\documentclass[10pt]{article}
\usepackage{latexsym, amsfonts, amsmath, amssymb, amscd, epsfig}

\numberwithin{equation}{section}
\newtheorem{thm}{Theorem}[section]

\newtheorem{re}{Remark}[section]
\newenvironment{pf}{{\noindent \it \bf Proof:}}{{\hfill$\Box$}\\}

\begin{document}

\title{A regularity criterion of 3D incompressible MHD system with mixed
pressure-velocity-magnetic field}
\author{Ahmad M. Alghamdi$^{1}$, Sadek Gala$^{2,3}$, Maria Alessandra Ragusa$%
^{3,4}$ \\
%EndAName
$^{{\small 1}}${\small Department of Mathematical Science , Faculty of
Applied Science,}\\
{\small Umm Alqura University, P. O. Box 14035, Makkah 21955, Saudi Arabia,}%
\\
$^{{\small 2}}${\small Department of Sciences Exactes, ENS of Mostaganem,
University of Mostaganem,}\\
\ {\small P.O.} {\small Box 227, Mostaganem 27000, Algeria,}\\
$^{{\small 3}}${\small Dipartimento di Matematica e Informatica, Universit%
\`{a} di Catania,}\\
{\small Viale Andrea Doria, 6 95125 Catania - Italy,}\\
$^{{\small 4}}${\small RUDN University, 6 Miklukho - Maklay St, Moscow,
117198, Russia}}
\date{}
\maketitle

\begin{abstract}
This work focuses on the 3D incompressible magnetohydrodynamic (MHD)
equations with mixed pressure-velocity-magnetic field in view of Lorentz
spaces. Our main result shows the weak solution is regular, provided that%
\begin{equation*}
{\frac{\pi }{\left( e^{-\left\vert x\right\vert ^{2}}+\left\vert
u\right\vert +\left\vert b\right\vert \right) ^{\theta }}\in L}%
^{p}(0,T;L^{q,\infty }(\mathbb{R}^{3})),\text{ \ where }\frac{2}{p}+\frac{3}{%
q}=2-\theta \text{ and }0\leq \theta \leq 1.
\end{equation*}
\end{abstract}

\noindent {\textit{Mathematics Subject Classification(2000):}\thinspace
\thinspace 35Q35,\thinspace \thinspace 35B65,\thinspace \thinspace 76D05}

\newpage

\section{Introduction}

\bigskip We are interested in the regularity of weak solutions to the
viscous incompressible magnetohydrodynamics (MHD) equations in $\mathbb{R}%
^{3}$
\begin{equation}
\left\{
\begin{array}{c}
\partial _{t}u+(u\cdot \nabla )u-\left( b\cdot \nabla \right) b-\Delta
u+\nabla \pi =0, \\
\partial _{t}b+(u\cdot \nabla )b-(b\cdot \nabla )u-\Delta b=0, \\
\nabla \cdot u=\nabla \cdot b=0, \\
u(x,0)=u_{0}(x),\text{ \ }b(x,0)=b_{0}(x),%
\end{array}%
\right.  \label{eq1.1}
\end{equation}%
where $u=(u_{1},u_{2},u_{3})$ is the velocity field, $b=(b_{1},b_{2},b_{3})$
is the magnetic field, and $\pi $ is the scalar pressure, while $u_{0}$ and $%
b_{0}$ are the corresponding initial data satisfying $\nabla \cdot
u_{0}=\nabla \cdot b_{0}=0$ in the sense of distribution.

Local existence and uniqueness theories of solutions to the MHD equations
have been studied by many mathematicians and physicists (see, e.g., \cite%
{CW, DL, ST}). But due to the presence of Navier-Stokes equations in the
system (\ref{eq1.1}) whether this unique local solution can exist globally
is an outstanding challenge problem. For this reason, there are many
regularity criteria of weak solutions for the MHD equations has been
investigated by many authors over past years (see e.g., \cite{DJZ, D, FJNZ,
G1, GR1, GR2, GRZ, LD, NGZ, Z1, Z2} and references therein). Note that the
literatures listed here are far from being complete, we refer the readers to
see for example \cite{GR20, JZ1, JZ2, JZ3, JZ4} for expositions and more
references.

More recently, Beir\~{a}o and Yang \cite{BY} proved the following regularity
criterion for the mixed pressure-velocity in Lorentz spaces for Leray-Hopf
weak solutions to 3D Navier-Stokes equations
\begin{equation}
{\frac{\pi }{\left( e^{-\left\vert x\right\vert ^{2}}+\left\vert
u\right\vert \right) ^{\theta }}\in L}^{p}(0,T;L^{q,\infty }(\mathbb{R}%
^{3})),\text{ \ where \ }0\leq \theta \leq 1\text{ and }\frac{2}{p}+\frac{3}{%
q}=2-\theta ,  \label{eq7}
\end{equation}%
where $L^{q,\infty }(\mathbb{R}^{3})$ denotes the Lorentz space (c.f. \cite%
{Tri}).

Motivated by the recent work of \cite{BY}, the purpose of this note is to
establish the regularity for the MHD equations (\ref{eq1.1}) with the mixed
pressure-velocity-magnetic in Lorentz spaces. Our main result can be stated
as follows:

\begin{thm}
\label{th1}Suppose that $(u_{0},b_{0})\in L^{2}(\mathbb{R}^{3})\cap L^{4}(%
\mathbb{R}^{3})$ with $\nabla \cdot u_{0}=\nabla \cdot b_{0}=0$ in the sense
of distribution.\ Let $\left( u,b\right) $ be a weak solution to the MHD
equations on some interval $\left[ 0,T\right] $ with $0<T\leq \infty $.\
Assume that $0\leq \theta \leq 1$ and that
\begin{equation}
{\frac{\pi }{\left( e^{-\left\vert x\right\vert ^{2}}+\left\vert
u\right\vert +\left\vert b\right\vert \right) ^{\theta }}\in L}%
^{p}(0,T;L^{q,\infty }(\mathbb{R}^{3})),\text{ \ where }\frac{2}{p}+\frac{3}{%
q}=2-\theta  \label{eq15}
\end{equation}%
then the weak $\left( u,b\right) $ is regular on $(0,T].$
\end{thm}

\begin{re}
A special consequence of Theorem \ref{th1} and its proof is the regularity
criterion of the 3D Navier-Stokes equations with the mixed pressure-velocity
in Lorentz spaces. This generalizes those of \cite{BY}.
\end{re}

In order to derive the regularity criterion of weak solutions to the MHD
equations (\ref{eq1.1}), we introduce the definition of weak solution.

Next, let us writing
\begin{equation*}
w^{\pm }=u\pm b,\ \ \text{\ }w_{0}^{\pm }=u_{0}\pm b_{0}.
\end{equation*}%
We reformulate equation (\ref{eq1.1}) as follows. Formally, if the first
equation of MHD equations (\ref{eq1.1}) plus and minus the second one,
respectively, then MHD equations (\ref{eq1.1}) can be re-written as:
\begin{equation}
\left\{
\begin{array}{l}
\partial _{t}w^{+}-\Delta w^{+}+(w^{-}\cdot \nabla )w^{+}+\nabla \pi =0, \\%
[3mm]
\partial _{t}w^{-}-\Delta w^{-}+(w^{+}\cdot \nabla )w^{-}+\nabla \pi =0, \\%
[2mm]
\mathrm{div}~w^{+}=0,~~~~\mathrm{div}~w^{-}=0, \\
w^{+}(x,0)=w_{0}^{+}(x),\text{ \ \ }w^{-}(x,0)=w_{0}^{-}(x).%
\end{array}%
\right.  \label{eq1.3}
\end{equation}%
The advantage is that the equations becomes symmetric.

\section{Proof of Theorem \protect\ref{th1}}

This section is devoted to the proof of Theorem \ref{th1}. In order to do
it, we first recall the following estimates for the pressure in terms of $u$
and $b$ (see e.g., \cite{GR20}) :
\begin{equation}
\left\Vert \pi \right\Vert _{L^{q}}\leq C\left( \left\Vert u\right\Vert
_{L^{2q}}^{2}+\left\Vert b\right\Vert _{L^{2q}}^{2}\right) ,\text{ \ \textrm{%
with} \ }1<q<\infty .  \label{eq120}
\end{equation}

We are now in position to prove our main result.

\begin{pf}
Multiplying the first and the second equations of (\ref{eq1.3}) by $%
\left\vert w^{+}\right\vert ^{2}w^{+}$ and $\left\vert w^{-}\right\vert
^{2}w^{-}$ , respectively, integrating by parts and summing up, we have
\begin{eqnarray*}
&&\frac{1}{4}\frac{d}{dt}(\left\Vert w^{+}\right\Vert
_{L^{4}}^{4}+\left\Vert w^{-}\right\Vert _{L^{4}}^{4})+\int_{\mathbb{R}%
^{3}}(\left\vert \nabla w^{+}\right\vert ^{2}\left\vert w^{+}\right\vert
^{2}+\left\vert \nabla w^{-}\right\vert ^{2}\left\vert w^{-}\right\vert
^{2})dx+\frac{1}{2}\int_{\mathbb{R}^{3}}(\left\vert \nabla \left\vert
w^{+}\right\vert ^{2}\right\vert ^{2}+\left\vert \nabla \left\vert
w^{-}\right\vert ^{2}\right\vert ^{2})dx \\
&=&-\int_{\mathbb{R}^{3}}\nabla \pi \cdot (w^{+}\left\vert w^{+}\right\vert
^{2}+w^{-}\left\vert w^{-}\right\vert ^{2})dx \\
&=&\int_{\mathbb{R}^{3}}\pi \cdot \mathrm{div}(w^{+}\left\vert
w^{+}\right\vert ^{2}+w^{-}\left\vert w^{-}\right\vert ^{2})dx \\
&\leq &\int_{\mathbb{R}^{3}}\left\vert \pi \right\vert (\left\vert
w^{+}\right\vert +\left\vert w^{-}\right\vert )(\nabla \left\vert
w^{+}\right\vert ^{2}+\nabla \left\vert w^{-}\right\vert ^{2})dx \\
&\leq &C\int_{\mathbb{R}^{3}}\left\vert \pi \right\vert ^{2}(\left\vert
w^{+}\right\vert +\left\vert w^{-}\right\vert )^{2}dx+\frac{1}{4}\int_{%
\mathbb{R}^{3}}(\left\vert \nabla \left\vert w^{+}\right\vert
^{2}\right\vert ^{2}+\left\vert \nabla \left\vert w^{-}\right\vert
^{2}\right\vert ^{2})dx.
\end{eqnarray*}%
Notice that $u=\frac{1}{2}(w^{+}+w^{-})$ and $b=\frac{1}{2}(w^{+}-w^{-})$,
then the above inequality means that
\begin{eqnarray}
&&\frac{d}{dt}(\left\Vert u\right\Vert _{L^{4}}^{4}+\left\Vert b\right\Vert
_{L^{4}}^{4})+2\left\Vert \nabla \left\vert u\right\vert ^{2}\right\Vert
_{L^{2}}^{2}+2\left\Vert \nabla \left\vert b\right\vert ^{2}\right\Vert
_{L^{2}}^{2}  \notag \\
&&+2\left\Vert \left\vert u\right\vert \left\vert \nabla u\right\vert
\right\Vert _{L^{2}}^{2}+2\left\Vert \left\vert b\right\vert \left\vert
\nabla b\right\vert \right\Vert _{L^{2}}^{2}+2\left\Vert \left\vert
u\right\vert \left\vert \nabla b\right\vert \right\Vert
_{L^{2}}^{2}+2\left\Vert \left\vert b\right\vert \left\vert \nabla
u\right\vert \right\Vert _{L^{2}}^{2}  \notag \\
&\leq &C\int_{\mathbb{R}^{3}}\left\vert \pi \right\vert ^{2}(\left\vert
u\right\vert +\left\vert b\right\vert )^{2}dx=K,  \label{eq21}
\end{eqnarray}%
where we have used
\begin{equation*}
\left\vert w^{+}\right\vert +\left\vert w^{-}\right\vert \leq \left\vert
w^{+}+w^{-}\right\vert +\left\vert w^{+}-w^{-}\right\vert .
\end{equation*}%
For $K$, borrowing the arguments in \cite{BY}, we set%
\begin{equation*}
V=e^{-\left\vert x\right\vert ^{2}}+\left\vert u\right\vert +\left\vert
b\right\vert \text{ \ \ and \ \ }\widetilde{\pi }={\frac{\pi }{\left(
e^{-\left\vert x\right\vert ^{2}}+\left\vert u\right\vert +\left\vert
b\right\vert \right) ^{\theta }}.}
\end{equation*}%
By the H\"{o}lder inequality and the following interpolation in Lorentz
space (see \cite{Tri})
\begin{equation*}
\left\Vert f^{\alpha }\right\Vert _{L^{p,q}(\mathbb{R}^{3})}\leq C\left\Vert
f\right\Vert _{L^{\alpha p,\alpha q}(\mathbb{R}^{3})}^{\alpha }\text{ \ \
for \ }\alpha >0,\text{ }p>0,\text{ }q>0,
\end{equation*}
we have%
\begin{eqnarray*}
K &=&\int_{\mathbb{R}^{3}}\left\vert \pi \right\vert ^{\lambda }V^{-\lambda
\theta }\left\vert \pi \right\vert ^{2-\lambda }V^{\lambda \theta
}(\left\vert u\right\vert +\left\vert b\right\vert )^{2}dx \\
&\leq &\int_{\mathbb{R}^{3}}\left\vert \widetilde{\pi }\right\vert ^{\lambda
}\left\vert \pi \right\vert ^{2-\lambda }V^{2+\lambda \theta }dx \\
&\leq &\left\Vert \left\vert \widetilde{\pi }\right\vert ^{\lambda
}\right\Vert _{L^{\frac{q}{\lambda },\infty }}\left\Vert \left\vert \pi
\right\vert ^{2-\lambda }\right\Vert _{L^{s,\frac{2}{2-\lambda }}}\left\Vert
V^{2\lambda }\right\Vert _{L^{r,\frac{2}{\lambda }}} \\
&=&\left\Vert \widetilde{\pi }\right\Vert _{L^{q,\infty }}^{\lambda
}\left\Vert \pi \right\Vert _{L^{s(2-\lambda ),2}}^{2-\lambda }\left\Vert
V^{2}\right\Vert _{L^{\lambda r,2}}^{\lambda },
\end{eqnarray*}%
where
\begin{equation*}
\frac{\lambda }{q}+\frac{1}{s}+\frac{1}{r}=1\text{ \ and \ }\lambda =\frac{2%
}{2-\theta }.
\end{equation*}%
By (\ref{eq120}), we have%
\begin{eqnarray*}
K &\leq &\left\Vert \widetilde{\pi }\right\Vert _{L^{q,\infty }}^{\lambda
}\left( \left\Vert \left\vert u\right\vert ^{2}\right\Vert _{L^{s(2-\lambda
),2}}+\left\Vert \left\vert b\right\vert ^{2}\right\Vert _{L^{s(2-\lambda
),2}}\right) ^{2-\lambda }\left\Vert V^{2}\right\Vert _{L^{\lambda
r,2}}^{\lambda } \\
&\leq &C\left\Vert \widetilde{\pi }\right\Vert _{L^{q,\infty }}^{\lambda
}\left\Vert V^{2}\right\Vert _{L^{s(2-\lambda ),2}}^{2-\lambda }\left\Vert
V^{2}\right\Vert _{L^{\lambda r,2}}^{\lambda }.
\end{eqnarray*}%
By the interpolation and Sobolev inequalities in Lorentz spaces, it follows
that%
\begin{equation}
\left\{
\begin{array}{c}
\left\Vert V^{2}\right\Vert _{L^{s(2-\lambda ),2}}\leq C\left\Vert
V^{2}\right\Vert _{L^{2,2}}^{1-\delta _{1}}\left\Vert V^{2}\right\Vert
_{L^{6,2}}^{\delta _{1}}\leq C\left\Vert V^{2}\right\Vert _{L^{2}}^{1-\delta
_{1}}\left\Vert \nabla V^{2}\right\Vert _{L^{2}}^{\delta _{1}}, \\
\left\Vert V^{2}\right\Vert _{L^{\lambda r,2}}\leq C\left\Vert
V^{2}\right\Vert _{L^{2,2}}^{1-\delta _{2}}\left\Vert V^{2}\right\Vert
_{L^{6,2}}^{\delta _{2}}\leq C\left\Vert V^{2}\right\Vert _{L^{2}}^{1-\delta
_{2}}\left\Vert \nabla V^{2}\right\Vert _{L^{2}}^{\delta _{2}},%
\end{array}%
\right.   \label{eq6.6}
\end{equation}%
where $0<\delta _{1},\delta _{2}<1$ and
\begin{equation*}
\frac{1}{s(2-\lambda )}=\frac{1-\delta _{1}}{2}+\frac{\delta _{1}}{6},\text{
\ }\frac{1}{\lambda r}=\frac{1-\delta _{2}}{2}+\frac{\delta _{2}}{6}.
\end{equation*}%
Hence from (\ref{eq6.6}) and Young inequality, it follows that%
\begin{eqnarray*}
K &\leq &C\left\Vert \widetilde{\pi }\right\Vert _{L^{q,\infty }}^{\lambda
}\left\Vert V^{2}\right\Vert _{L^{2}}^{(2-\lambda )(1-\delta _{1})+\lambda
(1-\delta _{2})}\left\Vert \nabla V^{2}\right\Vert _{L^{2}}^{(2-\lambda
)\delta _{1}+\lambda \delta _{2}} \\
&\leq &C\left\Vert \widetilde{\pi }\right\Vert _{L^{q,\infty }}^{\frac{%
2\lambda }{2-(2-\lambda )\delta _{1}-\lambda \delta _{2}}}\left\Vert
V^{2}\right\Vert _{L^{2}}^{2}+\frac{1}{2}\left\Vert \nabla V^{2}\right\Vert
_{L^{2}}^{2}.
\end{eqnarray*}%
Due to the definition of $V$, we see that
\begin{equation*}
\left\Vert V^{2}\right\Vert _{L^{2}}^{2}\leq C(1+\left\Vert \left\vert
u\right\vert +\left\vert b\right\vert \right\Vert _{L^{2}}^{2}+\left\Vert
\left\vert u\right\vert ^{2}+\left\vert b\right\vert ^{2}\right\Vert
_{L^{2}}^{2}),
\end{equation*}%
and%
\begin{equation*}
\left\Vert \nabla V^{2}\right\Vert _{L^{2}}^{2}\leq C(1+\left\Vert
\left\vert u\right\vert +\left\vert b\right\vert \right\Vert
_{L^{2}}^{2}+\left\Vert \nabla (\left\vert u\right\vert +\left\vert
b\right\vert )\right\Vert _{L^{2}}^{2}+\left\Vert \nabla (\left\vert
u\right\vert ^{2}+\left\vert b\right\vert ^{2})\right\Vert _{L^{2}}^{2}).
\end{equation*}%
Consequently, we get%
\begin{eqnarray*}
K &\leq &C\left\Vert \widetilde{\pi }\right\Vert _{L^{q,\infty }}^{\frac{%
2\lambda }{2-(2-\lambda )\delta _{1}-\lambda \delta _{2}}}(1+\left\Vert
\left\vert u\right\vert +\left\vert b\right\vert \right\Vert
_{L^{2}}^{2}+\left\Vert \left\vert u\right\vert ^{2}+\left\vert b\right\vert
^{2}\right\Vert _{L^{2}}^{2}) \\
&&+C(1+\left\Vert \left\vert u\right\vert +\left\vert b\right\vert
\right\Vert _{L^{2}}^{2}+\left\Vert \nabla (\left\vert u\right\vert
+\left\vert b\right\vert )\right\Vert _{L^{2}}^{2})+\frac{1}{2}\left\Vert
\nabla (\left\vert u\right\vert ^{2}+\left\vert b\right\vert
^{2})\right\Vert _{L^{2}}^{2} \\
&\leq &C\left\Vert \widetilde{\pi }\right\Vert _{L^{q,\infty }}^{\frac{%
2\lambda }{2-(2-\lambda )\delta _{1}-\lambda \delta _{2}}}(1+\left\Vert
u\right\Vert _{L^{2}}^{2}+\left\Vert b\right\Vert _{L^{2}}^{2}+\left\Vert
u\right\Vert _{L^{4}}^{4}+\left\Vert b\right\Vert _{L^{4}}^{4}) \\
&&+C(1+\left\Vert u\right\Vert _{L^{2}}^{2}+\left\Vert b\right\Vert
_{L^{2}}^{2}+\left\Vert \nabla u\right\Vert _{L^{2}}^{2}+\left\Vert \nabla
b\right\Vert _{L^{2}}^{2})+\frac{1}{2}\left\Vert \nabla \left\vert
u\right\vert ^{2}\right\Vert _{L^{2}}^{2}+\frac{1}{2}\left\Vert \nabla
\left\vert b\right\vert ^{2}\right\Vert _{L^{2}}^{2}.
\end{eqnarray*}%
Since $(u,b)$ is a weak solution to (\ref{eq1.1}), then $(u,b)$ satisfies%
\begin{equation*}
(u,b)\in L^{\infty }(0,T;L^{2}(\mathbb{R}^{3}))\cap L^{2}(0,T;H^{1}(\mathbb{R%
}^{3})).
\end{equation*}%
Inserting the above estimates into (\ref{eq21}), we obtain
\begin{eqnarray*}
&&\frac{d}{dt}(\left\Vert u\right\Vert _{L^{4}}^{4}+\left\Vert b\right\Vert
_{L^{4}}^{4})+\left\Vert \nabla \left\vert u\right\vert ^{2}\right\Vert
_{L^{2}}^{2}+\left\Vert \nabla \left\vert b\right\vert ^{2}\right\Vert
_{L^{2}}^{2} \\
&&+2\left\Vert \left\vert u\right\vert \left\vert \nabla u\right\vert
\right\Vert _{L^{2}}^{2}+2\left\Vert \left\vert b\right\vert \left\vert
\nabla b\right\vert \right\Vert _{L^{2}}^{2}+2\left\Vert \left\vert
u\right\vert \left\vert \nabla b\right\vert \right\Vert
_{L^{2}}^{2}+2\left\Vert \left\vert b\right\vert \left\vert \nabla
u\right\vert \right\Vert _{L^{2}}^{2} \\
&\leq &C\left\Vert \widetilde{\pi }\right\Vert _{L^{q,\infty }}^{\frac{%
2\lambda }{2-(2-\lambda )\delta _{1}-\lambda \delta _{2}}}(1+\left\Vert
u\right\Vert _{L^{2}}^{2}+\left\Vert b\right\Vert _{L^{2}}^{2}+\left\Vert
u\right\Vert _{L^{4}}^{4}+\left\Vert b\right\Vert _{L^{4}}^{4}) \\
&&+C(1+\left\Vert u\right\Vert _{L^{2}}^{2}+\left\Vert b\right\Vert
_{L^{2}}^{2}+\left\Vert \nabla u\right\Vert _{L^{2}}^{2}+\left\Vert \nabla
b\right\Vert _{L^{2}}^{2}) \\
&\leq &C\left\Vert \widetilde{\pi }\right\Vert _{L^{q,\infty }}^{\frac{%
2\lambda }{2-(2-\lambda )\delta _{1}-\lambda \delta _{2}}}(1+\left\Vert
u\right\Vert _{L^{4}}^{4}+\left\Vert b\right\Vert
_{L^{4}}^{4})+C(1+\left\Vert \nabla u\right\Vert _{L^{2}}^{2}+\left\Vert
\nabla b\right\Vert _{L^{2}}^{2}),
\end{eqnarray*}%
Using Gronwall's inequality with the assumption (\ref{eq15}), we deduce that%
\begin{equation*}
(u,b)\in L^{\infty }(0,T;L^{4}(\mathbb{R}^{3}))\subset L^{8}(0,T;L^{4}(%
\mathbb{R}^{3})).
\end{equation*}%
We complete the proof of Theorem \ref{th1}.
\end{pf}

\end{document}